\documentclass{ifacconf}

\usepackage{graphicx}      % include this line if your document contains figures
\usepackage{natbib}        % required for bibliography

\usepackage{amsmath,amssymb,amsfonts,bm}
\usepackage{bbm}
\usepackage{subfigure}
\usepackage{textcomp}
\usepackage{algorithm}
\usepackage{algpseudocode}

\allowdisplaybreaks[4]
\usepackage{url}
\begin{document}
\begin{frontmatter}

\title{Globally Convergent Policy Gradient Methods for Linear Quadratic Control of Partially Observed Systems} 

\thanks[footnoteinfo]{The research was supported by National Natural Science Foundation of China under Grant no. 62033006, and Tsinghua University Initiative Scientific Research Program.}

\author{Feiran Zhao}, 
\author{Xingyun Fu},  
\author{Keyou You} 
%\author[Third]{Third C. Author}

\address{Department of Automation and BNRist, Tsinghua University, Beijing 100084, China (e-mail: zhaofr18@mails.tsinghua.edu.cn; fxy20@mails.tsinghua.edu.cn; youky@tsinghua.edu.cn.)}

\begin{abstract}           
While the optimization landscape of policy gradient methods has been recently investigated for {partially observed} linear systems in terms of both static output feedback and dynamical controllers, they only provide convergence guarantees to stationary points. In this paper, we propose a new policy parameterization for {partially observed} linear systems, using a past input-output trajectory of finite length as feedback. We show that the solution set to the parameterized optimization problem is a matrix space, which is invariant to {similarity transformation}. By proving a gradient dominance property, we show the global convergence of policy gradient methods. Moreover, we observe that the gradient is orthogonal to the solution set, revealing an explicit relation between the resulting solution and the initial policy. Finally, we perform simulations to validate our theoretical results.
\end{abstract}


\end{frontmatter}
%===============================================================================

%% There are a number of predefined theorem-like environments in
%% ifacconf.cls:
%%
%% \begin{thm} ... \end{thm}            % Theorem
%% \begin{lem} ... \end{lem}            % Lemma
%% \begin{claim} ... \end{claim}        % Claim
%% \begin{conj} ... \end{conj}          % Conjecture
%% \begin{cor} ... \end{cor}            % Corollary
%% \begin{fact} ... \end{fact}          % Fact
%% \begin{hypo} ... \end{hypo}          % Hypothesis
%% \begin{prop} ... \end{prop}          % Proposition
%% \begin{crit} ... \end{crit}          % Criterion

\section{Introduction}

Recent years have witnessed tremendous successes of reinforcement learning (RL) in applications such as sequential decision-making  problems~\citep{mnih2015human-level, silver2016mastering} and continuous control \citep{tobin2017domain, levine2016end, andrychowicz2020learning, recht2019tour}. As an essential approach of RL, the policy gradient (PG) method directly searches over a policy space to optimize a performance index of interests using sampled trajectories without any identification process. Such an end-to-end approach is conceptually simple and easy to implement in practice.

In contrast to the above empirical successes, the theoretical understanding of the PG method has largely lagged as it often involves challenging non-convex optimization problems. To fill this gap, there has been a resurgent interest in studying the theoretical properties of PG methods for classical control problems \citep{fazel2018global,gravell2020learn,zhao2022sample,zhao2021global, malik2019derivative, zhang2021policy,li2019distributed,zheng2021analysis,zheng2022escaping,duan2022optimization,duan2021optimization,fatkhullin2021optimizing}.
The seminal work of \cite{fazel2018global} has shown that the well-known linear quadratic regulator (LQR) problem~\citep{zhou1996robust} has a gradient dominance property, leading to the global convergence of PG methods despite the non-convexity. There have also been other PG-based works considering, e.g., system stabilization~\citep{zhao2022sample}, robustness~\citep{zhang2021policy} and distributed control~\citep{li2019distributed}, just to name a few. We refer the readers to the survey~\citep{bin2022towards} for a comprehensive overview.

In this paper, we consider \textit{partially observed} linear systems, where the state cannot be directly observed and only input-output trajectories are available as feedback. For partially observed systems, different PG methods have been investigated in~\citep{zheng2021analysis,zheng2022escaping,duan2022optimization,duan2021optimization,fatkhullin2021optimizing}. Depending on how the control policy is parameterized, they can be broadly categorized into static output feedback (SOF)~\citep{fatkhullin2021optimizing, duan2021optimization} and dynamic output feedback~\citep{zheng2021analysis,zheng2022escaping,duan2022optimization}. The former class only uses the current output as feedback, while the latter uses all past input-output trajectories by invoking a linear filter. In both classes, the optimization landscape of PG methods can be substantially different from that in state feedback control. Particularly, the gradient dominance property does not hold, which is the key to the convergence in~\cite{fazel2018global}. Moreover, the set of stabilizing controllers is usually disconnected, and stationary points can be local minima or saddle points~\citep{fatkhullin2021optimizing, duan2021optimization}. Even though a perturbed PG method is proposed to escape the strict saddle, its convergence rate has not been well characterized yet~\citep{zheng2021analysis,zheng2022escaping}. Last but not least, the cost function may vary with similarity transformations~\citep{duan2022optimization}, which further increases the difficulty in the convergence analysis. Therefore, all the above PG methods for partially observed linear systems can only provide convergence guarantees to stationary points.

In this paper, we propose a PG method for linear quadratic control with global convergence for {partially observed} linear systems. We first propose a new policy parameterization in the form of input-output feedback (IOF) using a \textit{past} input-output trajectory of fixed length instead of the \textit{current} output. Then, we show that the solution set to the parameterized optimization problem is a matrix space, which is invariant to the similarity transformation. Even though an optimal policy is not unique, our problem still meets the gradient dominance property, based on which we prove the global convergence of the PG method. Moreover, we reveal an explicit relation between the solution and the initial policy by observing that the gradient is orthogonal to the solution set. We also propose a zero-order algorithm with warm-up cost evaluation for sample-based implementation. Finally, we perform simulations to validate our theoretical results.

The remainder of this paper is organized as follows. Section \ref{sec:formu} formulates the linear quadratic control problem for partially observed systems as a parameterized optimization problem. Section \ref{sec:optimal} derives its optimal solution which is shown to be invariant to similarity transformation. Section \ref{sec:land} shows the convergence of the PG method and discusses its implementation in the sample-based setting. Section \ref{sec:sim} presents a numerical case study. Conclusion is made in Section \ref{sec:conc}.

\section{Problem formulation} \label{sec:formu}
Consider the partially observed linear system
\begin{equation}\label{equ:sys}
\begin{aligned}
x_{t+1} &= Ax_t + Bu_t, \\
y_t &= Cx_t,
\end{aligned}
\end{equation}
where $x_t \in \mathbb{R}^n$ is the state, $u_t \in \mathbb{R}^m$ is the control input, and $y_t \in \mathbb{R}^d$ is the measurable output. The matrices $(A,B,C)$ are model parameters.

We aim to find a policy sequence $\{\pi_t\}$ using only past input-output data to minimize an infinite-horizon quadratic cost, i.e.,
\begin{equation}\label{prob:lqg}
\begin{aligned}
\mathop{\text { minimize }}\limits_{\{\pi_t\}}&~J(\{\pi_t\}):= \mathbb{E}_{x_0 \sim \mathcal{D}} \left[\sum_{t=0}^{\infty}(y_{t}^{\top} Q y_{t}+u_{t}^{\top} R u_{t})\right]\\
\text {subject to} & ~(\ref{equ:sys}), ~ u_t = \pi_t(u_{-\infty},y_{-\infty},\dots, u_{t-1}, y_{t-1})
\end{aligned}
\end{equation}
with $Q\geq 0, R>0$. We require the distribution $\mathcal{D}$ of the initial state $x_0$ to satisfy the following assumption.

\begin{assum}\label{assum:distri}
	 The distribution $\mathcal{D}$ has zero mean with a positive definite covariance matrix $\Sigma_0 = \mathbb{E}[x_0x_0^{\top}]>0$.
\end{assum}

Since the control policy does not depend on the current output $y_t$, it can be applied to strictly causal systems. Throughout the paper, we make the following assumption standard in the control theory~\citep{zhou1996robust}. 
\begin{assum}
	\label{assumption}
	$(A,B)$ are controllable and $(C, A)$ are observable.
\end{assum}

When the state is measurable, the optimal policy to (\ref{prob:lqg}) is linear state feedback 
\begin{equation}\label{equ:state_feedback}
u_t = - (R+B^{\top}P^*B)^{-1}B^{\top}P^*A x_t, 
\end{equation}
where $P^*$ is the unique positive semi-definite solution to the algebraic Riccati equation (ARE)~\citep{zhou1996robust}
$$
P^*=A^{\top} P^* A+Q_c-A^{\top} P^* B(R+B^{\top} P^* B)^{-1} B^{\top} P^* A
$$
with $Q_c = C^{\top} Q C$. Note that the optimal policy (\ref{equ:state_feedback}) is independent of $\Sigma_0$. 

This paper considers a new input-output feedback (IOF) policy parameterization
\begin{equation}\label{equ:policy}
u_t = -Kz_{t,p},
\end{equation}
where ${z}_{t,p} = [{u}_{t,p}^{\top}, {y}_{t,p}^{\top}]^{\top}$, ${u}_{t,p} = [u_{t-1}^{\top}, \cdots, u_{t-p}^{\top}]^{\top}$, ${y}_{t,p} = [y_{t-1}^{\top}, \cdots, y_{t-p}^{\top}]^{\top}$, $p \in \mathbb{N}$ is a system-dependent constant to be defined later, and $K \in \mathbb{R}^{m\times q}$ with $q = p(m+d)$ is the gain matrix. The intuition behind (\ref{equ:policy}) is that the state can be recovered from a finite-length past input-output trajectory under Assumption \ref{assumption}. 

In this paper, we use gradient methods to solve the following problem viewing $K$ as the optimization matrix
\begin{equation}\label{prob:opt}
\mathop{\text { minimize }}\limits_{K}~J(K),~~
\text {subject to}~K \in \mathcal{S},
\end{equation}
where $J(K)$ is the quadratic cost following the policy (\ref{equ:policy}) and $\mathcal{S}$ is the feasible set containing all the stabilizing policy. Clearly, this is a challenging constrained non-convex optimization problem. In fact, an optimal solution to (\ref{prob:opt}) may not be unique, which makes (\ref{prob:opt}) more challenging to solve. In the sequel, we investigate the optimization landscape of (\ref{prob:opt}) to show the global convergence of policy gradient methods.

%We show that the solution set to (\ref{prob:opt}) is a matrix space, which is invariant to similarity transformations. This key feature renders its optimization landscape substantially distinct from that of the LQG~\citep{duan2022optimization,zheng2021analysis} and SOF~\citep{fatkhullin2021optimizing,duan2021optimization} problems.

% Interestingly, we show that the gradient is orthogonal to the solution set, which along with the gradient dominance property leads to the global convergence of PG methods.

%This key feature enables us to prove a gradient dominance property that is however missing in \cite{zheng2021analysis,duan2021optimization,duan2022optimization,fatkhullin2021optimizing} and further show the global convergence of PG methods.

\section{Optimal IOF control policy}\label{sec:optimal} 
This section shows that the solution set to (\ref{prob:opt}) is a matrix space, which is invariant to similarity transformation.

We first express the state $x_t$ using the trajectory $z_{t,p}$. Let $o$ and $c$ be the observability index and controllability index, respectively, and $p = \max\{o,c\}\leq n$. Then, the following matrices
\begin{equation}\label{equ:obse}
\mathcal{O}_p = \begin{bmatrix}
CA^{p-1} \\
\vdots \\
CA \\
C
\end{bmatrix}, ~\text{and}~\mathcal{C}_p = [B~~AB~~\cdots~~A^{p-1}B]
\end{equation}
have full column and row rank, respectively. At time step $t$, the state can be represented using system dynamics and history trajectories as
\begin{equation}\label{equ:traj}
\begin{aligned}
x_t &= A^p x_{t-p} + \mathcal{C}_p{u}_{t,p} \\
{y}_{t,p} &= \mathcal{O}_p x_{t-p} + \mathcal{T}_p{u}_{t,p}
\end{aligned}
\end{equation}
with a Toeplitz matrix
$$
\mathcal{T}_p = \begin{bmatrix}
0 & C B & C A B & \cdots & C A^{p-2} B \\
0 & 0 & C B & \cdots & C A^{p-3} B \\
\vdots & \vdots & \ddots & \ddots & \vdots \\
0 & \cdots & & 0 & C B \\
0 & 0 & 0 & 0 & 0
\end{bmatrix}.
$$
In some cases, $p$ is unknown, and we only have knowledge of the system order $n$, input dimension $m$, and output dimension $d$. Then, one can substitute $p$ with $n$ in (\ref{equ:obse}), and the results in this paper still hold. For simplicity, we omit the subscript $p$ where it can be understood from the context. 

Since $\mathcal{O}$ has full column rank, it has a unique left pseudo inverse $\mathcal{O}^{\dagger} = (\mathcal{O}^{\top}\mathcal{O})^{-1}\mathcal{O}^{\top} $. Then, it follows immediately from (\ref{equ:traj}) that $x_t$ can be uniquely determined by eliminating $x_{t-p}$ as
\begin{equation}\label{equ:xt}
x_t = (\mathcal{C} - A^p\mathcal{O}^{\dagger}\mathcal{T}) {u}_{t,p} + A^p\mathcal{O}^{\dagger}{y}_{t,p} := S {z}_{t,p}
\end{equation}
with $S = [\mathcal{C} - A^pO^{\dagger}\mathcal{T}, A^pO^{\dagger}]$. Clearly, $S$ has full row rank by noting
\begin{equation}\notag
S \begin{bmatrix}
I & 0\\
\mathcal{T} & I
\end{bmatrix} =  [\mathcal{C} - A^pO^{\dagger}\mathcal{T}, A^pO^{\dagger}]\begin{bmatrix}
I & 0\\
\mathcal{T} & I
\end{bmatrix} = [\mathcal{C}, A^pO^{\dagger}],
\end{equation}
and has a unique right pseudo inverse $S^{\dagger} = S^{\top}(SS^{\top})^{-1}$.

Then, the feasible set of $K$ can be written as
$$
	\mathcal{S} = \{K \in\mathbb{R}^{m\times q}| \rho(A-BKS^{\dagger}) <1 \},
$$
where $\rho(\cdot)$ denotes the spectral radius of a square matrix. We have the closed-form expression for $J(K)$.
%Suppose that (\ref{equ:xt}) holds for all time step $t = 0,1,\dots$ (may require that $\{u_{-p}, y_{-p}, \dots, u_{-1}, y_{-1}\}$ is available, to be discussed later), 
\begin{lem}
	For any $K \in \mathcal{S}$,
	the cost function $J(K)$ can be written as 
	$$
	J(K) = \text{Tr}(P_K\Sigma_0),
	$$
	where $P_K\geq 0$ is the solution to the Lyapunov equation
	\begin{equation}\label{equ:std_lya}
	\begin{aligned}
	P_K = &Q_c + (S^{\dagger})^{\top}K^{\top}RKS^{\dagger} \\
	&+ (A - BKS^{\dagger})^{\top}P_K(A-BKS^{\dagger}).
	\end{aligned}
	\end{equation}	
\end{lem}
\begin{pf}
Let $V_K(x) = x^{\top}P_Kx$  be the value function of problem (\ref{prob:opt}) following the stabilizing policy $K$. By the well-known Bellman equation~\citep{bertsekas1995dynamic}, it follows that
$$
V_K(x_t) = y_t^{\top}Qy_t + (-Kz_{t,p})^{\top}R(-Kz_{t,p}) + V_K(x_{t+1}).
$$

Then, substituting $x_t$ with (\ref{equ:xt}) yields that
\begin{align*}
z_{t,p}^{\top}S^{\top}P_KSz_{t,p} = z_{t,p}^{\top}S^{\top}Q_cSz_{t,p} + z_{t,p}^{\top}K^{\top}RKz_{t,p}\\ + z_{t,p}^{\top}(AS - BK)^{\top}P_K(AS-BK)z_{t,p}.
\end{align*}

Noting that it holds for all $z_{t,p}$, it holds that
$$
	S^{\top}P_KS = S^{\top}Q_cS + K^{\top}RK + (AS - BK)^{\top}P_K(AS-BK).
$$

Pre- and post-multiplying $(S^{\dagger})^{\top}$ and $S^{\dagger}$ in both sides of the above equation yields (\ref{equ:std_lya}). Since by definition $J(K)= \mathbb{E}_{x_0}[V_K(x_0)]$, we complete the proof. \qed
\end{pf}

Clearly, an optimal policy of form (\ref{equ:policy}) can be determined by substituting (\ref{equ:xt}) into (\ref{equ:state_feedback}) as
$
u_t = - K^*{z}_{t,p},
$
with 
\begin{equation}\label{equ:opt_K}
K^* = (R+B^{\top}P^*B)^{-1} B^{\top}P^*A S,
\end{equation}
which also satisfies the following Lyapunov equation
\begin{equation}\label{equ:Lya_opt}
\begin{aligned}
P^* = &Q_c + (S^{\dagger})^{\top}(K^*)^{\top}RK^*S^{\dagger} \\
&+ (A - BK^*S^{\dagger})^{\top}P^*(A-BK^*S^{\dagger}).
\end{aligned}
\end{equation}
However, an optimal solution is not unique as $S^{\dagger}$ does not have full row rank. Define the matrix space 
$$\mathcal{F} = \{\Delta\in\mathbb{R}^{m\times q} | \Delta \cdot S^{\dagger}=0 \}.$$
We show in the following theorem that the solution set to (\ref{prob:opt}) is a matrix space parallel to $\mathcal{F}$.
\begin{thm}\label{theorem:optimal}
	Define the set $\mathcal{K} = \{K \in \mathbb{R}^{m\times q} | K = K^* + \Delta, \Delta \in \mathcal{F} \}$. Then, $K$ is an optimal policy to (\ref{prob:opt}) if and only if $K \in \mathcal{K}$.
\end{thm}
\begin{pf}
	To prove the ``if" statement, let $K \in \mathcal{K}$. By the definition of $\mathcal{K}$, it holds that $KS^{\dagger} = (K^* + \Delta)S^{\dagger} = K^*S^{\dagger}, \forall \Delta \in \mathcal{F}$. Combining with (\ref{equ:std_lya}), all $K \in \mathcal{K}$ have the same optimal cost as $K^*$, which implies that $K \in \mathcal{K}$ is optimal.
	
	To prove the ``only if" statement, suppose that $K$ is optimal, i.e., $K$ satisfies
	$$
	P^* = Q_c + (S^{\dagger})^{\top}K^{\top}RKS^{\dagger} + (A - BKS^{\dagger})^{\top}P^*(A-BKS^{\dagger}).
	$$
	
	Let $KS^{\dagger} = K^*S^{\dagger} + E$. Taking $KS^{\dagger}$ into the above equation leads to that
	\begin{align*}
		&(K^*S^{\dagger})^{\top}RE + E^{\top}RK^*S^{\dagger} + E^{\top}(R+B^{\top}P^*B)E \\
		&= (A-BK^*S^{\dagger})^{\top}P^*BE  + (BE)^{\top}P^*(A-BK^*S^{\dagger}).
	\end{align*}
	
	Then, inserting (\ref{equ:opt_K}) into the above equation yields that
	$E^{\top}(R+B^{\top}P^*B)E = 0$. Thus, we can only have $E = 0$, i.e., $K \in \mathcal{K}$. The proof is now completed. \qed
\end{pf}

At the first sight, Theorem \ref{theorem:optimal} appears to be a negative result, as the global convergence of PG methods has been shown in the existing literature only when an optimal policy is unique. In fact, the convergence in our problem can be proved by utilizing the special structure of the solution set, as to be shown in Section \ref{sec:land}.

%which can also be solved via the Bellman optimality equation
%$$
%x_t^{\top}P^*x_t = \min_{u_t} \{x_t^{\top}Q_cx_t + u_t^{\top}Ru_t + x_{t+1}^{\top}P^*x_{t+1}\}.
%$$

Finally, we show that the optimal solution set is invariant to similarity transformations.

\begin{lem}\label{lem:similarity}
	For a nonsingular matrix $T$, define the new system with the similarity transformation $\widetilde{x}_t = Tx_t$
	\begin{equation}\label{equ:new_sys}
	\begin{aligned}
	\widetilde{x}_{t+1} &= TAT^{-1}\widetilde{x}_t + TBu_t, \\
	y_t &= CT^{-1}\widetilde{x}_t.
	\end{aligned}
	\end{equation}	
 Then, $K$ is an optimal policy of (\ref{equ:new_sys}) if and only if $K \in \mathcal{K}$.
\end{lem}
\begin{pf}
	The ARE of (\ref{equ:new_sys}) can be written as
	\begin{align*}
	&0=(TAT^{-1})^{\top} \widetilde{P}^* TAT^{-1}-\widetilde{P}^*+(CT^{-1})^{\top}QCT^{-1}\\
	&-(TAT^{-1})^{\top} \widetilde{P}^* TB(R+(TB)^{\top} \widetilde{P}^* TB)^{-1} (TB)^{\top}\\
	&\times  \widetilde{P}^* TAT^{-1}.
	\end{align*}
	Pre- and post-multiplying $T^{\top}$ and $T$, it follows that $\widetilde{P}^* = (T^{-1})^{\top}P^*T^{\top}$. Similarly, we can show that $\widetilde{S} = TS$.
	
	Hence, an optimal gain matrix of the new system is 
	\begin{align*}
	\widetilde{K}^* &= (R+(TB)^{\top}\widetilde{P}^*TB)^{-1} (TB)^{\top}\widetilde{P}^*TAT^{-1} \widetilde{S} \\
	&= (R+B^{\top}P^*B)^{-1} B^{\top}P^*A S\\
	&= K^*.
	\end{align*}
	Noting $\mathcal{N}((\widetilde{S}^{\dagger})^{\top}) = \mathcal{N}((S^{\dagger})^{\top})$, the proof is completed. \qed
\end{pf}

The optimal dynamical control policy in \cite{duan2022optimization, zheng2021analysis} has the following form 
\begin{equation}
\begin{aligned}
&\dot{\xi}=(A-B K) \xi+L(y-C \xi) \\
&u=-K \xi,
\end{aligned}
\end{equation}
where $K$ is the LQR gain, $L$ is the Kalman gain~\citep{zhou1996robust}, and $\xi$ is the internal state. Clearly, it is not unique and each similarity transformation to (\ref{equ:sys}) leads to a different optimal policy, which makes it much more challenging to provide any convergence guarantees. In contrast, Lemma \ref{lem:similarity} implies that we can focus on the minimal realization in (\ref{equ:sys}) to study the optimization landscape of the PG method.

\section{Policy gradient method for IOF control}\label{sec:land}
In this section, we propose a PG method under IOF parameterization with global convergence. Then, we present a zero-order optimization algorithm to solve an optimal policy by only using sampled trajectories.

\subsection{Global convergence}
Define $\Sigma_K$ as the solution to the Lyapunov equation
$$
\Sigma_K = \Sigma_0 + (A-BKS^{\dagger})\Sigma_K(A-BKS^{\dagger})^{\top}.
$$
Then, we have the following gradient expression.
\begin{figure}[t]
	\centerline{\includegraphics[width=60mm]{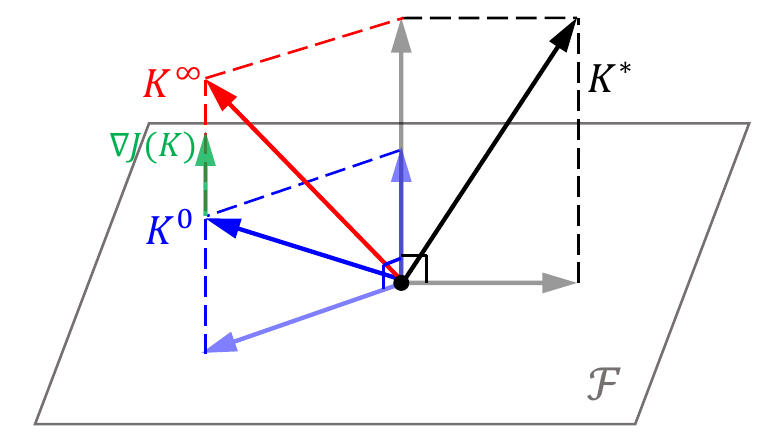}}
	\caption{Subspace relations among $K^0, K^{\infty}, K^*$ and $\nabla J(K)$.}
	\label{pic:opt_policy}
\end{figure}
\begin{lem}
	For $K \in \mathcal{S}$, the gradient of $J(K)$ is 
	$$
	\nabla J(K) = 2E_K\Sigma_K(S^{\dagger})^{\top},
	$$
	where $E_{K}=(R+B^{\top} P_K B) KS^{\dagger}-B^{\top} P_K A$. 
\end{lem}
\begin{pf}
	Let $X = KS^{\dagger}$. By \citet[Theorem 1]{fazel2018global}, the gradient with respect to $X$ can be written as
	$$
	\nabla_{X} J = 2E_K\Sigma_K.
	$$
	
	Then, it follows from the chain rule that
	$$
	\nabla J(K) = \nabla_{X} J \cdot  (S^{\dagger})^{\top}=2E_K\Sigma_K(S^{\dagger})^{\top}. \qed
	$$ 
\end{pf}

Consider the following gradient method to update $K$
\begin{equation}\label{equ:pg}
K^{i+1} = K^i - \eta \nabla J(K^i),~ i \in \{0,1,\dots\}.
\end{equation}
As in the standard LQR~\citep{fazel2018global}, we show that $J(K)$ has the following gradient dominance property (aka Polyak-Lojasiewicz condition~\citep{polyak1963gradient}), which guarantees that all stationary points are optimal. 

\begin{lem}\label{lem:gd_domi}
For any $K \in \mathcal{S}$, it holds that
$$
J(K)-J(K^*)\leq \frac{\|\Sigma^*\|\|S\|^2}{4\underline\sigma(R)\underline\sigma^2(\Sigma_K)}\operatorname{tr} \{
\nabla{J}(K)^{\top}\nabla{J}(K)	
\},
$$
where $\Sigma^*$ denotes $\Sigma_{K^*}$ and $\underline{\sigma}(\cdot)$ denotes the smallest eigenvalue of a square matrix.
\end{lem}
\begin{pf}
	By \citet[Corollary 5]{fazel2018global}, we can show that the cost satisfies
	\begin{align*}
	&J(K)-J(K^*)\leq \frac{\|\Sigma^*\|}{4\underline\sigma(R)\underline\sigma^2(\Sigma_K)}\operatorname{tr} \{
	\nabla_X{J}^{\top}\nabla_X{J}	
	\} \\
	&= \frac{\|\Sigma^*\|}{4\underline\sigma(R)\underline\sigma^2(\Sigma_K)}\operatorname{tr} \{S\nabla{J}(K)^{\top}\nabla{J}(K)S^{\top}\}\\
	&\leq \frac{\|\Sigma^*\|\|S\|^2}{4\underline\sigma(R)\underline\sigma^2(\Sigma_K)}\operatorname{tr} \{
	\nabla{J}(K)^{\top}\nabla{J}(K)	
	\}. \qed
	\end{align*}
\end{pf}

It can be observed that the gradient is orthogonal to the matrix space $\mathcal{F}$. Define $\Pi_{\mathcal{F}}$ as the projection operator of a matrix onto $\mathcal{F}$ and $\Pi_{\mathcal{F}}^{\perp}$ onto its orthogonal space. 

\begin{lem}\label{lem:ortho}
	Let $K \in \mathcal{S}$. Then, we have $\Pi_{\mathcal{F}}(\nabla J(K)) = 0$.
\end{lem}
\begin{pf}
For any $\Delta \in \mathcal{F}$, it holds that
$$
\text{tr}\{\Delta^{\top}\cdot \nabla J(K) \} = 2\text{tr}\{ E_K\Sigma_K(S^{\dagger})^{\top}\Delta^{\top}\}=0.
$$
Hence, $\nabla J(K)$ is orthogonal to $\mathcal{F}$. \qed
\end{pf}
This fact implies that for any initial policy $K^0\in\mathcal{S}$, its projection $\Pi_{\mathcal{F}}(K^0)$ will not be affected by the update (\ref{equ:pg}). Along with Lemma \ref{lem:gd_domi}, we have the following global convergence guarantees.

\begin{thm}\label{theorem:pg}
	For $K^0\in \mathcal{S}$ and an appropriate stepsize $\eta$ that is polynomial in problem parameters, e.g., $\|A\|$, $\|B\|$, $\|S\|$, ${\underline{\sigma}(\Sigma_0)}$, $\underline{\sigma}(Q)$, $\underline{\sigma}(R)$, the gradient update (\ref{equ:pg})
	converges to  
	\begin{equation}\label{equ:conv_policy}
			K^{\infty} := \lim\limits_{i\rightarrow \infty} K^i = \Pi_{\mathcal{F}}(K^0) + \Pi_{\mathcal{F}}^{\perp}(K^*)
	\end{equation}
	at a linear rate, i.e., for $i \in \{0,1,\dots\}$,
	$$
	J(K^{i+1}) - J^* \leq \left(1 - \frac{ 2\eta \sigma_0^2 \underline{\sigma}(R)}{\|\Sigma^*\|\|S\|^2}\right)(J(K^i) - J^*).
	$$
\end{thm}
\begin{pf}
	The convergence follows the same vein as the proof of \cite[Theorem 7]{fazel2018global} based on Lemma \ref{lem:gd_domi}, and is omitted here due to space limitation. By Lemma \ref{lem:ortho} and Theorem \ref{theorem:optimal}, it follows that
	$
	\Pi_{\mathcal{F}}(K^{\infty}) = \Pi_{\mathcal{F}}(K^0)$ and $
	\Pi_{\mathcal{F}}^{\perp}(K^{\infty}) = \Pi_{\mathcal{F}}^{\perp}(K^*),
	$
	leading to (\ref{equ:conv_policy}). \qed
\end{pf}

Fig. \ref{pic:opt_policy} illustrates the subspace relations among $K^0, K^{\infty}, K^*$ and $\nabla J(K)$. Theorem \ref{theorem:pg} ensures that for any initial stabilizing policy $K^0 \in \mathcal{S}$, the PG update in (\ref{equ:pg}) converges to the solution set $\mathcal{K}$ at a linear rate. Our convergence rate depends on $\|S\|$, which tends to infinity as the smallest eigenvalue of the observability matrix $\mathcal{O}$ tends to zero. This implies that the system needs to be ``sufficiently observable".

In contrast to the static output feedback (SOF) parameterization \citep{duan2022optimization, duan2021optimization}, we do not require the observation matrix $C$ to have full row rank~\citep{duan2022optimization} or full column rank~\citep{duan2021optimization} to prove the global convergence. Even though both optimal IOF and dynamical control policies in \cite{duan2022optimization, zheng2021analysis} are not unique, our PG method has global convergence due to the gradient dominance property and the invariance to similarity transformation.

In the following, we discuss the implementation of our PG method when an explicit model $(A,B,C)$ is unavailable.

\subsection{Sample-based implementation}\label{sec:implement}

In the sample-based setting, the gradient can only be estimated via zero-order information. However, it is challenging to evaluate the cost function, as implementing $u_0 = -Kz_{0,p}$ requires $\{u_{-p},y_{-p},\dots, u_{-1},y_{-1}\}$ to be known.

To generate the required sequence, we use a random control policy $u_t \sim \mathcal{N}(0,I)$ in $t \in \{-p,\dots,-1\}$. More specifically, we generate a trajectory by
\begin{equation}\label{equ:sample}
\begin{aligned}
&x_{-p}\sim \mathcal{N}(0,I)  ~\text{and}~ \\
&u_t = \begin{cases} w_t, w_t\sim \mathcal{N}(0,I), & t\in \{-p,\dots,-1\} \\ -Kz_{t,p}, & t \in \{0,\dots,T\} \end{cases}
\end{aligned}
\end{equation}
where $\{w_t\}_{-p}^{-1}$ is an independent random sequence. By the dynamics in (\ref{equ:sys}), $x_0$ satisfies
$$
x_0 = \begin{bmatrix}
\mathcal{C}_p & A^p
\end{bmatrix} \begin{bmatrix}
u_{0,p}\\
x_{-p}
\end{bmatrix}.
$$
Since $\mathcal{C}_p$ has full row rank, the distribution of $x_0$ generated by (\ref{equ:sample}) satisfies Assumption \ref{assum:distri}.

Then, we estimate the cost by a single sampled trajectory
\begin{equation}\label{equ:warm}
\widehat{J}(K) = \sum_{t=0}^{T}(y_{t}^{\top} Q y_{t}+u_{t}^{\top} R u_{t})\\
\end{equation}
following (\ref{equ:sample}). In practice, the sampling time $T$ can be set sufficiently large to well approximate the infinite-horizon cost. We refer to (\ref{equ:warm}) as warm-up cost evaluation since it involves another (random) policy in $t = \{-p,\dots,-1\}$. 

We present our zero-order algorithm with warm-up cost evaluation in Algorithm \ref{alg:learning}. Particularly, we use a two-point method~\citep{malik2019derivative} to estimate the gradient in step 2-5. The parameter $r$ is called the smoothing radius used to control the variance of the gradient estimate. For the convergence analysis of Algorithm \ref{alg:learning}, one can apply standard results~\citep[Theorem 1]{malik2019derivative} of zero-order methods, and is omitted in this paper.
\begin{algorithm}[t]
	\caption{The zero-order algorithm with two-point gradient estimate}
	\label{alg:learning}
	\begin{algorithmic}[1]
		\Require
		An initial policy $K^0 \in \mathcal{S}$, the number of iterations $N$, a smoothing radius $r$, the stepsize $\eta$.
		\For{$i=0,1,\cdots , N-1$}
		\State Sample a perturbation matrix $U^i$ uniformly from the unit sphere $\mathcal{U}^{mq-1}$.
		\State Set ${K}_1^i = K^{i}+ r\sqrt{mq}U^{i}$ and ${K}_2^i = K^{i} - r\sqrt{mq}U^{i}$.
		\State Obtain $\widehat{J}({K}_1^i)$ and $\widehat{J}({K}_2^i)$ from (\ref{equ:warm}).
		\State Estimate the gradient $$\widehat{\nabla J^i} = \frac{1}{2r}(\widehat{J}({K}_1^i)-\widehat{J}({K}_2^i))U^{i}.$$
		\State Update the policy by $K^{i+1} = K^{i} - \eta \widehat{\nabla J^i}$.
		\EndFor	
		\Ensure A policy $K^{N}$.
	\end{algorithmic}
\end{algorithm}
\section{Simulation}\label{sec:sim}
In this section, we validate the convergence of our PG methods via simulations. Moreover, we compare the performance of our IOF control policy from Algorithm \ref{alg:learning} with the SOF policy. The simulation code is provided in \url{https://github.com/fuxy16/Input-output-Feedback}.
\begin{figure}[t]
	\centerline{\includegraphics[width=70mm]{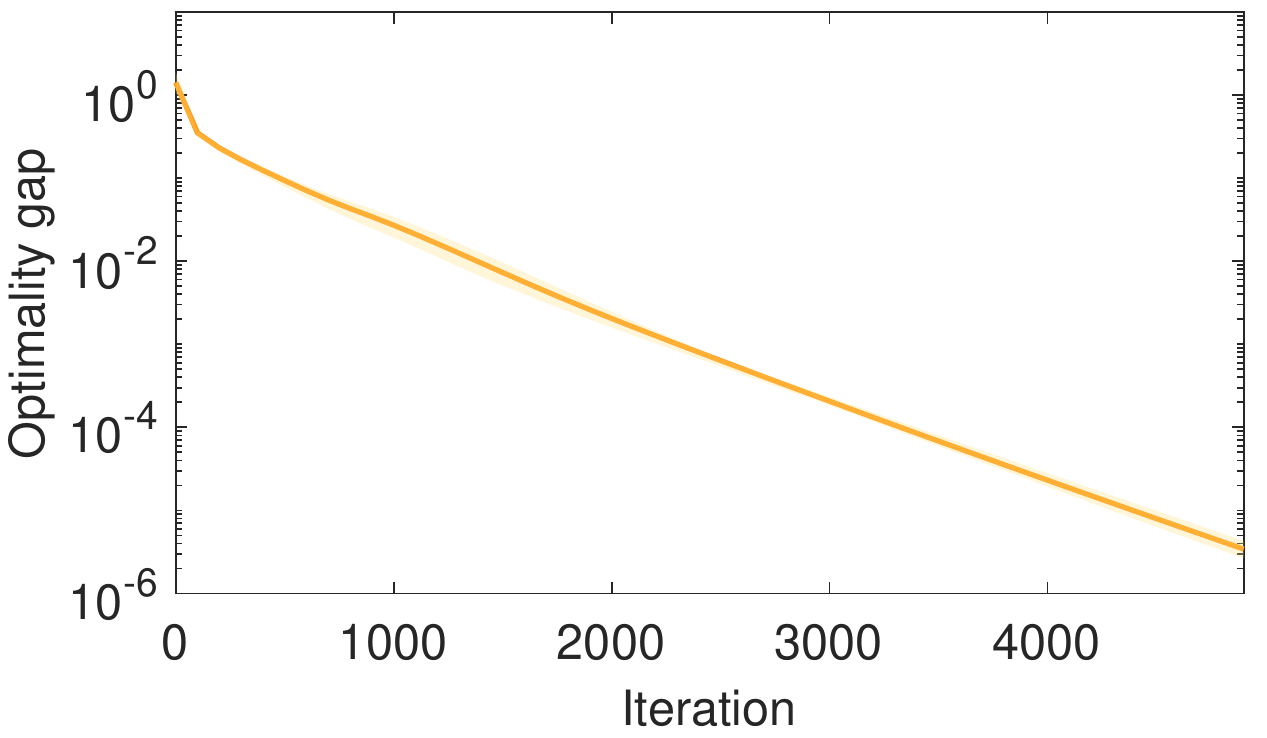}}
	\caption{Convergence of the PG update in \eqref{equ:pg}.}
	\label{pic:mb}
\end{figure}
\subsection{Example}
To validate the convergence, we randomly generate a dynamical model $(A,B,C)$ with $n=4,m=d=2$ as
$$
A = \begin{bmatrix}
0.568  &  0.215  &  0.122  & -0.156\\
-0.074  & -0.021 &  -0.114 &  -0.307\\
0.568  &  0.211  &  0.047  & -0.604\\
-0.455  &  1.141 &  -0.204 &  -0.478
\end{bmatrix}
$$
\begin{equation}\notag
\begin{aligned}
&B = \begin{bmatrix}
0.584  &  1.193\\
-0.988 &   0.696\\
0.176  & -0.683\\
0.470  & -1.163
\end{bmatrix}\\
&C = \begin{bmatrix}
0.719  & -0.138   & 1.026 &  -0.743\\
-1.014  &  0.252 &  -0.500 &  -0.601
\end{bmatrix}.
\end{aligned}
\end{equation}
This open-loop stable system is both controllable and observable with $p=2$. Let $Q = I_4$ and $R = 0.01\times I_2$. Then, Assumption \ref{assumption} is satisfied. In the sequel, we search over the matrix space $\mathbb{R}^{2\times8}$ to find an optimal solution. 

\subsection{Convergence of our PG method}
In the model-based setting, we perform (\ref{equ:pg}) using the model $(A,B,C)$ to validate the global convergence result in Theorem \ref{theorem:pg}. Let the stepsize be $\eta = 10^{-3}$. The initial policy $K^0$ is selected by first generating a random matrix with its elements being Gaussian and then normalizing it such that $\rho(A-BK^0)=0.8$. Then, we conduct (\ref{equ:pg}) where the gradient is computed using $(A,B,C)$ and display the optimality gap in the cost $(J(K^i)-J^*)/J^*$ in Fig. \ref{pic:mb}. The bold centreline denotes the mean of 20 independent trials and the shaded region demonstrates their variance. As expected from Theorem \ref{theorem:pg}, the gap diminishes fast at a linear rate, and the randomness of $K^0$ only induces a small variance.
\begin{figure}[t]
	\centerline{\includegraphics[width=70mm]{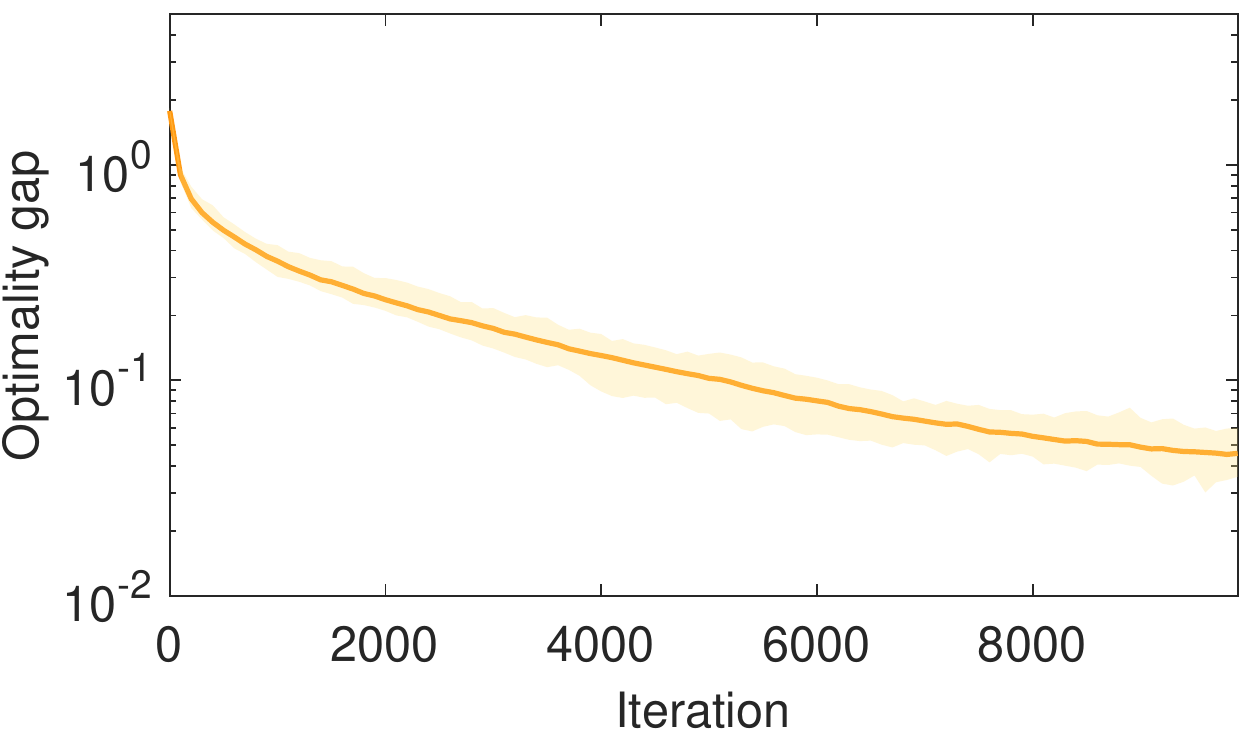}}
	\caption{Convergence of Algorithm \ref{alg:learning}.}
	\label{pic:mf}
\end{figure}
\begin{table}[!t]\label{table1}
	\caption{Average cost of the resulted policy.}
	\begin{center}
		\begin{tabular}{|cccc|}
			\hline
			& Optimal &IOF & SOF   \\
			\hline
			Model-based (d=2) & $8.282$ & $\bf{8.987}$ & $13.013$\\
			\hline
			Sample-based (d=2) & $10.502$ &$\bf{11.438}$ & $15.519$ \\
			\hline
			Model-based (d=4) & $19.649$ & $\bf{19.793}$ & $20.611$\\
			\hline
			Sample-based (d=4) & $26.112$ &$\bf{26.428}$ & $29.467$\\
			\hline
		\end{tabular}
	\end{center}
\end{table}
%
%\begin{table}[!t]\label{table2}
%	\caption{The region of attraction radius $r_{\text{roa}}$.}
%	\begin{center}
%		\begin{tabular}{|cccc|}
%			\hline
%			& Optimal &Ours & SOF   \\
%			\hline
%			Model-based & $19.649$ & $\bf{19.793}$ & $20.611$\\
%			\hline
%			Sample-based & $26.112$ &$\bf{26.428}$ & $29.467$ \\
%			\hline
%		\end{tabular}
%	\end{center}
%\end{table}

In the sample-based setting, we perform Algorithm \ref{alg:learning} to demonstrate the performance of our zero-order method. Set the sampling time $T = 20$, the stepsize $\eta = 10^{-5}$, the smoothing radius $r = 0.2$ and $K^0 = 0$. The convergence is shown in Fig. \ref{pic:mf}, where the variance originates from the warm-up process of (\ref{equ:warm}).

\subsection{Comparison with SOF policy}
To show the merits of our new parameterization, we compare it with the SOF control policy in the form of
\begin{equation}
	u_t = -K_sy_t,
\end{equation}
where $K_s$ is solved by PG methods in \cite{duan2021optimization}. Particularly, we consider two cases, $d=2$ and $d=4$. When $d=2$, the matrix $C$ is rank deficient and $K_s$ is only guaranteed to be locally minimal~\citep{duan2021optimization, polyak1963gradient}. When $d=4$, we sample a new model $(A,B,C)$ randomly where the matrix $C$ is invertible. Hence, the state $x_t$ can be recovered by $x_t = C^{-1}y_t$ and $K_s$ is expected to have the same performance as the LQR control~\citep[Theorem 1]{duan2021optimization}. We set the stepsize by grid search which is $\eta = 10^{-5}$ for both IOF and SOF, and set other parameters as before. We perform $10^5$ iterations of PG updates and compare the performance between the resulting IOF and SOF policies. Their average of infinite-horizon costs in 20 independent trials are displayed in Table 1. In the case $d=2$, the results are reasonable as the SOF only converges to local minima. Surprisingly, even in the case $d=4$ the IOF policy still yields a lower cost. We note that the matrix space of the SOF problem is $\mathbb{R}^{2\times 4}$, which is $\mathbb{R}^{2\times 12}$ for the IOF problem. Thus, this result means that our PG method converges faster even though its gain matrix has a higher dimension.

\section{Conclusion}\label{sec:conc}
In this paper, we have proposed a new parameterization of the policy for partially observed linear systems, under which the PG method has been shown to globally converge to the solution set. We have also found some interesting properties such as the orthogonality of the gradient, and the invariance of the solution set to similarity transformation.

We now discuss some possible future works. Since this paper only considers the vanilla gradient descent method, it would be interesting to investigate the performance of both natural gradient and Gauss-Newton methods, which have been shown to have a faster convergence rate in the LQR problem \citep{fazel2018global}. It is also interesting to see whether the convergence can be preserved in the presence of process and measurement noises. 
 
\bibliography{mybibfile}             % bib file to produce the bibliography
                                                     % with bibtex (preferred)
%\appendix
%\section{A summary of Latin grammar}    % Each appendix must have a short title.
%\section{Some Latin vocabulary}              % Sections and subsections are supported  
                                                                    
\end{document}